\documentclass[oneside,reqno]{amsart}
\usepackage{amsthm}
\usepackage{amssymb}
\usepackage{mathrsfs}
\usepackage{mathtools}
\usepackage{url}		
\usepackage{comment}

\newtheorem{theorem}{Theorem}[section]
\newtheorem{lemma}[theorem]{Lemma}
\newtheorem{conjecture}{Conjecture}

\usepackage{comment}
\numberwithin{equation}{section}

\begin{document}

\title[Twin primes]{ On the distribution of twin primes}
\author{Madieyna Diouf}
\address{Madieyna Diouf \\ Arizona State University,
Tempe, AZ 85282}
\email{mdiouf1@asu.edu}
\thanks{}



\date{08/08/2021}


\begin{abstract}
We introduce a sieve for counting twin primes up to a given range. Our method 
depends on a parameter ${\lambda}_x$ and the estimation of the number of twin primes obtained as a result, is called a fundamental structure of the distribution of twin primes. Combining the latter with an asymptotic bound of ${\lambda}_x$, establishes venues, conducive to a discovery of a partial result that can be considered as a suitable variant of the prime number theorem. Furthermore, we obtain an asymptotic bound of the number of twin primes less than $x$.  
\end{abstract}
\maketitle

       \section{Introduction And Statement of Results}
The following statements describe our findings including key quantitative and technical results on a sieve for counting twin primes. A fundamental structure of the distribution of twin primes is  obtained by sieve in $(2.8)$
\begin{align}
{\pi}_2(x) \approx \frac{\pi(x)({\lambda}_x+1)-x}{{\lambda}_x-1}
\end{align}
where ${{\lambda}_x}$ is the average number of composite integers between the pairs of consecutive primes whose gap exceeds two and its lower member is less than $x$.\\
This opens possibilities to reach some asymptotic and non-asymptotic bounds of the number of primes and twin primes up to $x$. Consequently, we prove the following variant of the prime number theorem by combining respectively 
\begin{align}
\: \: \: \: \: \: \: \: \: \: \: \: \text{Theorem $4.2$\: }& &\pi(x) \sim \frac{x}{{{\lambda}_x}},\: \: \: \: \: \: \: \: \: \: \: \: \: \: \: \: \: \: \: \: \: \: \: \: \: \: \: \: \: \: \: \: \: \: \\
\: \: \: \: \: \: \: \: \: \: \: \: \text{and Theorem $5.1$\: }& &{\lambda}_x \sim \log x. \: \: \: \: \: \:\: \: \: \: \: \: \: \: \: \: \: \: \: \: \: \: \: \: \: \: \: \: \: \: \: \: \: 
\end{align}
The aforementioned results are explicit and unconditional. By using $(1.1)$ and invoking $(1.3)$, one has the main result in Theorem $3.1$  
\begin{align}
{\pi}_2(x) \sim (2\pm o(\log x))\frac{x}{{\log}^2 x}.
\end{align}
Clearly $(1.4)$ establishes an asymptotic bound of ${\pi}_2(x)$. Our method revolves around a sieve that counts twin primes by enumerating their midpoints.
\subsection{Background} 
Let's first give a brief history of the twin prime conjecture and outline some recent breakthroughs making dramatic new progress on questions related to the conjecture. 
\begin{conjecture}[Twin Prime Conjecture] 
There are infinitely many twin primes. 
\end{conjecture}
G. H. Hardy and J.E. Littlewood \cite{HarLit} formulated a similar, but stronger twin prime conjecture known as the first Hardy-Littlewood conjecture. 
\begin{conjecture}[Strong Twin Prime or first Hardy-Littlewood conjecture] 
\begin{align*}
{\pi}_2(x) \sim 2 C_2 \frac{x}{(\log x)^2} &\sim 2{C}_2\int_{2}^x \frac{dt}{{\log}^2t}.\\
\text{Here ${C}_2$ refers to the twin prime }&\text{constant.}\\
{C}_2= \prod_{\textstyle{p\;{\rm prime}\atop p \ge 3}}&\left(1 - \frac{1}{(p-1)^2}\right)= 0.66016 ...
\end{align*}
\end{conjecture} 
In $1849$, the Twin Prime conjecture was considered by the French Mathematician Alphonse de Polignac \cite{Pol},  but there has been speculation that it could go back as Euclid and the ancient Greeks over 2000 years ago. Attempts to prove the conjecture have produced partial results that stand on their own merits. Results including Brun's theorem \cite{Brun}, Mertens' theorems \cite{Mer} and the Hardy-Littlewood conjecture \cite{Har}, along with Chen's theorem \cite{Chen}.\\ 
Currently, the result obtained by sieve on twin primes is an upper bound that states the following. 
There exists a constant $K>0$ such that for $x\geq 3$
\begin{align*}
{\pi}_2(x) &\leq 2K{C}_2\int_{2}^x \frac{dt}{{\log}^2t}.
\end{align*}
The best value of $K$ obtained in this direction is slightly below $3.4$ whenever $x$
is large enough (see \cite{Wu}, building on the earlier work of \cite{Che}, \cite{Fou}, \cite{Fouu}, \cite{Bom}, \cite{Fouv}).\\
The original work of the pioneer of modern Sieve Theory Viggo Brun,  implied a convergence of the sum of reciprocals of twin primes. \\
Chen has proved that There are infinitely many primes $p$ such that $p+2$ has at most $2$ prime factors. Unfortunately, {on a march} toward the twin prime conjecture, his method faces the parity problem that prevents sieves from giving good estimates.\\
Bounded gaps between primes is, in recent years, the most used alternative approximation to the Twin Prime Conjecture. 
In the groundbreaking work of Goldston-Pintz-Yıldırım \cite{Gol}, 
they proved that given $\epsilon>0$, there are infinitely many pairs $(p_1, p_2)$ of distinct primes such that 
\begin{align*}
|p_2-p_1|\leq \epsilon \: \log p_1.
\end{align*}
Their method also showed that if the Elliott-Halberstam conjecture is true, then the gap above can be reduced to $16$. Despite many significant results in this direction \cite{EFou}, \cite{Bom}, \cite{Bom2},  \cite{Bom3}, it is not until $2013$ that Zhang \cite{Zhang}  proved an acceptable alternative of this conjecture, which allowed him to use the GPY method to prove that there are infinitely many bounded gaps between primes. 
 Zhang's work validates the existence of infinitely many pairs $(p_1, p_2)$ of distinct primes such that 
\begin{align*}
|p_2-p_1|\leq 70 000 000.
\end{align*}
In less than a year after Zhang's results, a new variant of the GPY method was
discovered independently by Maynard \cite{May} and Tao (unpublished).\\

A great deal of effort has already been expended on the Twin Prime conjecture. 
It is widely believed that new ideas may be required to prove the conjecture. \\ A fair question is,  how does the idea presented here differ from those of the previous authors? The improvement comes from a new sieve for finding twin primes. It offers a platform, where with a little bit of effort, one can make a significant contribution to the twin prime conjecture by engineering good bounds of ${\lambda}_x$ for $(1.1)$. 
\section{A sieve for counting twin primes}
\subsection{Preliminary Steps}
Consider a composite odd integer $x$; fixed, once for all. \\
$\bullet$ Let $p_n$ be the largest prime not exceeding $x$. \\
$\bullet$ We denote by $\pi_2(x)$, the number of pairs of consecutive primes whose gap is two or less. These are the twin primes and the pair $(2, 3)$.\\
$\bullet$ Say $s_i$ is a sequence of composite integers between two consecutive primes whose gap exceeds two and its lower member is less than $x$. \\
Set $j$ to be the number of these sequences $s_i$ less than $x$. Thus,
\begin{center}
$s_1=(8, 9, 10)$, $s_2=(14,15, 16)$, $s_3=(20, 21, 22)$, ..., $s_j=(p_n+1,..., x)$.
\end{center}
By the definition of $s_i$, the last sequence denoted by $s_j$, is expected to be \\$s_j=(p_n+1, ..., p_{n+1}-1)$. But since $p_n<x<p_{n+1}$ and we are only interested in the integers $\leq x$, then the last sequence shall end at $x$ to be $s_j=(p_n+1, ..., x)$. \\
$\bullet$ Let $\mathscr{S}_1$ be the set of sequences $s_i$, for $i$ going from $1$ to $j$,
\begin{center}
$\mathscr{S}_1=\{(8, 9, 10), (14,15, 16), (20, 21, 22), ..., (p_n+1, ..., x)\}.$
\end{center}
$\bullet$ We define $\mathscr{S}_2$ as the set of pairs of consecutive primes whose gap exceeds two and its lower member is less than $x$,
\begin{center}
$\mathscr{S}_2=\{_{_{7}}( )_{_{11}}\: ,\:  _{_{13}}( )_{_{17}}\: , \: _{_{19}}( )_{_{23}},\: ..., \: _{_{p_n}}()_{_{p_{n+1}}} \}.$
\end{center}
The notation $_{_{p_i}}()_{_{p_{i+1}}}$ indicates a pair of consecutive primes $p_i$ and $p_{i+1}$.\\ 
$\bullet$ Observe that $\mathscr{S}_1$ and $\mathscr{S}_2$ are different sets with the same cardinality. That is
\begin{align*}
|\mathscr{S}_1|=|\mathscr{S}_2|. 
\end{align*}
$\bullet$ Let ${{\lambda}_x}$ hold the average number of composites per sequence $s_i$ in $\mathscr{S}_1$. Clearly, this is the same definition of ${\lambda}_x$ as given in the introduction. Meaning that ${{\lambda}_x}$ points to the average number of composites between the pairs of consecutive primes whose gap exceeds two and its lower member is less than $x$. 
\subsection{Sieving for twin primes $\leq x$}
\text{}\\
At the presence of the positive integers less than or equal to $x$,
\begin{center}
$1, 2, 3, 4, 5, 6, 7, 8, 9, 10, 11, 12, 13, 14, 15, 16, 17, 18, 19, 20, 21, 22, 23, 24, ..., x.$
\end{center}
\textbf{a)} Remove all primes $\leq x$.
\begin{center}
$1,\: \: ,\: \: ,\: 4,\: \: ,\: 6,\: \: ,\: 8, 9, 10,\: \: \: , 12,\: \: \: , 14, 15, 16,\: \: \: , 18,\: \: \: , 20, 21, 22,\:...,\: p_n+1,....,\: x$
\end{center}
\textbf{b)} Remove all elements in the $\cup_{s_i\in \mathscr{S}_1}$.
\begin{center}
$1,\: \: ,\: \: ,\: 4,\: \: ,\: 6,\: \: ,\: \: ,\: \: , \: \: ,\: \: \: , 12,\: \: \: , \: \: , \: \: , \: \: ,\: \: \: , 18,\: \: \: , \: \: , \: \: , \: \: ,\: \: \: , \: \: , \: \: .... ... $
\end{center}
The remaining integers are the midpoints of the twin primes $<x$, except the positive integer $1$. We wish to count these midpoints since each of them represents a unique pair of twin primes. The positive integer $1$ shall count for the midpoint of the first pair $(2, 3)$ which is not a twin prime and its midpoint is not represented in the sieve. Proceedings of steps $a)$ and $b)$ are expressed by:
\begin{align}
x-\pi(x)-|{\cup}_{{s_i}\in\mathscr{S}_1}| &={\pi}_2(x).\nonumber\\
|{\cup}_{{s_i}\in\mathscr{S}_1}|&=x-\pi(x)-{\pi}_2(x).
\end{align}
Where
\begin{align}
 |{\cup}_{{s_i}\in\mathscr{S}_1}|&\approx {{\lambda}_x} \cdot|\mathscr{S}_1|.
\end{align}
The approximation symbol $\lq\lq \approx "$ in $(2.2)$ and in the remainder of this paper, shall be defined as \lq\lq almost equal to".  It is used above to indicate that the real value of $|{\cup}_{{s_i}\in\mathscr{S}_1}|$ is almost equal to its estimation ${{\lambda}_x} \cdot|\mathscr{S}_1|$. In $(2.2)$, we state that the number of integers in the union of sequences $s_i$ elements of the set $\mathscr{S}_1$, is approximately equal to the number of sequences in the set $\mathscr{S}_1$ that is $|\mathscr{S}_1|$, multiplied by the average number of integers per sequence denoted by ${\lambda}_x$.\\ 
With $(2.1)$ and $(2.2)$ we get
\begin{align}
{\lambda}_x \cdot|\mathscr{S}_1|&\approx x-\pi(x)-{\pi}_2(x).\nonumber\\
|\mathscr{S}_1|&\approx \frac{x-\pi(x)-{\pi}_2(x)}{{{\lambda}_x}}.
\end{align}
\textbf{c)} Given the primes $< x$.
\begin{center}
$2, 3, 5, 7, 11, 13, 17, 19, 23, ..., p_n$
\end{center}
Remove the (lower members) of the pairs of consecutive primes whose gap exceeds two. 
\begin{center}
$\: \: \: \: 2, 3, 5, \: \: , 11, \: \: , 17, \: \: , \: \: ....$
\end{center}
The remaining must be (the lower members) of the pairs whose gap is two or less; This gives the value of ${\pi}_2(x)$ by definition. Step $c)$ is formulated as
\begin{align}
\pi(x)-|\mathscr{S}_2|&=\pi_2(x).
\end{align}
Since 
\begin{align}
|\mathscr{S}_2|=|\mathscr{S}_1|,
\end{align}
then substituting $|\mathscr{S}_2|$ in $(2.4)$ with $|\mathscr{S}_1|$ gives
\begin{align}
|\mathscr{S}_1|&=\pi(x)-\pi_2(x).
\end{align}
With $(2.3)$ and $(2.6)$, we obtain
\begin{align}
 \pi(x)-\pi_2(x) &\approx \frac{x-\pi(x)-{\pi}_2(x)}{{\lambda}_x}.\\
{\pi}_2(x)\big({\lambda}_x-1\big)&\approx \pi(x)\big({\lambda}_x+1\big)-x.\nonumber\\
{\pi}_2(x)&\approx \frac{\pi(x)\big({\lambda}_x+1\big)-x}{{\lambda}_x-1}.
\end{align}
\subsection{Quick Example}
Say $x=25$, then
\begin{align}
\pi(x)&=9.\\
\log(x)&=3.22.\\
x/\log x&=7.77.\\
{\lambda}_x &=|{\cup}_{{s_i}\in\mathscr{S}_1}|/|\mathscr{S}_1|.\\
&=11/4=2.75.\\
x/{\lambda}_x&=25/2.75=9.09.
\end{align}
The value of $x/{\lambda}_x$ in this example is a better approximation of $\pi(x)$ than $x/\log x$. In addition, the number of twin primes $<25$ (including the pair $(2,3)$) is ${\pi}_2(x)=5$.
When using estimations obtained by sieve, we have
\begin{align}
{\pi}_2(x)&\approx \frac{\pi(x)\big({\lambda}_x+1\big)-x}{{\lambda}_x-1}.\\
&\approx \frac{9\big(2.75+1\big)-25}{2.75-1}.\\
{\pi}_2(x)&\approx 5.
\end{align}
The estimation of ${\pi}_2(x)$ in $(2.8)$ frames,  a fundamental structure of the distribution of twin primes. It suffices now to find a good bound for ${\lambda}_x$. 
\section{An asymptotic bound of ${\lambda}_x$}
This section is devoted to proving that ${\lambda}_x \sim \log x$. It requires Lemma $3.1$ and $3.2$ together with Theorem $3.3$, before the result is exhibited in Theorem $3.4$.
\begin{align}
\text{Set \: \: \: \: \: \: \: \: \: \: \: \: \: \: \: \: \: \: } k_x=\pi(x)/{\pi}_2(x).\: \: \: \: \: \: \: \: \: \: \: \: \: \: \: \: \: \: \: \: \: \: \: \: \: \: \: \: \: \: \: \: \: \: \: \: \: \: \: \: \: \: \: \: \: \: \: \: \: \: \: \: \: \: 
\end{align}
Let $c_x$ be the real value that satisfies the equation
\begin{align}
\pi(x)=\frac{x}{\log x-c_x} \text{\: \:  (for the sake of precision).}
\end{align}
It is well known that $c_x \to 1$ when $x \to \infty$ as the result of the prime number theorem. \\Going forward, the symbol $\lq\lq \lesssim "$ shall be read as \lq\lq less than or almost equal to".
\begin{lemma}
\[\frac{k_x}{k_x-1}\log x-7 \lesssim \: {\lambda}_x \: \lesssim \frac{k_x}{k_x-1}\log x-1.\]
\end{lemma}
\begin{proof}
By following the steps in the sieve and looking for ${\lambda}_x$ instead of ${\pi}_2(x)$, \\we obtain by $(2.7)$
\begin{align}
{\lambda}_x &\approx \frac{x-\pi(x)-{\pi}_2(x)}{\pi(x)-{\pi}_2(x)}.\\
{\lambda}_x &\approx \frac{x}{\pi(x)-{\pi}_2(x)}-\frac{\pi(x)}{\pi(x)-{\pi}_2(x)}-\frac{{\pi}_2(x)}{\pi(x)-{\pi}_2(x)}.\\
{\lambda}_x &\approx \frac{x}{\pi(x)-{\pi}(x)/k_x}-\frac{\pi(x)}{\pi(x)-{\pi}(x)/k_x}-\frac{{\pi}(x)/k_x}{\pi(x)-{\pi}(x)/k_x}.\\
{\lambda}_x &\approx \frac{k_x x}{(k_x-1)\pi(x)}-\frac{k_x}{k_x-1}-\frac{1}{k_x-1}.
\end{align}
In $(3.6)$, replace $\pi(x)$ with its expression given in $(3.2)$, we get
\begin{align}
{\lambda}_x &\approx \frac{k_x (\log x-c_x)}{k_x-1}-\frac{k_x}{k_x-1}-\frac{1}{k_x-1}.\\
{\lambda}_x &\lesssim \frac{k_x (\log x-c_x)}{k_x-1}-\frac{k_x}{k_x-1}.
\end{align}
By $(3.7)$
\begin{align}
{\lambda}_x &\gtrsim  \frac{k_x (\log x-c_x)}{k_x-1}-\frac{k_x}{k_x-1}-1.
\end{align}
With $(3.8)$ and $(3.9)$ we have 
\begin{align}
\frac{k_x (\log x-c_x)}{k_x-1}-\frac{k_x}{k_x-1}-1 \lesssim \: &{\lambda}_x \: \lesssim \frac{k_x (\log x-c_x)}{k_x-1}-\frac{k_x}{k_x-1}.\\
\frac{k_x}{k_x-1}(\log x-c_x-1)-1 \lesssim \: &{\lambda}_x \: \lesssim \frac{k_x}{k_x-1}(\log x-c_x-1).\\
\frac{k_x}{k_x-1}(\log x-1)-\frac{k_xc_x}{k_x-1}-1 \lesssim \: &{\lambda}_x \: \lesssim \frac{k_x}{k_x-1}(\log x-1),
\end{align}
where
\begin{align}
\frac{k_x}{k_x-1}<2.
\end{align}
By Dusart \cite{Dus}: If $x>60184$, then 
\begin{align}
\frac {x} {\log x - 1} &< \pi(x) <  \frac {x} {\log x - 1.1}.
\end{align}
A combination of $(3.14)$ and $(3.2)$ implies that
\begin{align}
1< c_x<1.1 \text{\: \: \: for $x>60184$}.
\end{align}
With $(3.13)$ and $(3.15)$, one has
\begin{align}
\frac{k_xc_x}{k_x-1}<4.
\end{align}
By using $(3.16)$ to simplify the left inequality in $(3.12)$, we get
\begin{align}
\frac{k_x}{k_x-1}(\log x-1)-5 \lesssim \: &{\lambda}_x \: \lesssim \frac{k_x}{k_x-1}(\log x-1).
\end{align}
So that
\begin{align}
\frac{k_x}{k_x-1}\log x-7 \lesssim \: &{\lambda}_x \: \lesssim \frac{k_x}{k_x-1}\log x-1.
\end{align}
\end{proof}
\begin{lemma}
\[\frac{x}{{\lambda}_x+2} \lesssim \pi(x)\: \lesssim \frac{k_xx}{{\lambda}_x(k_x-1)+k_x}.\]
\end{lemma}
\begin{proof}
Of Lemma $3.1$, we deduce via $(3.11)$ that
\begin{align}
\frac{k_x}{k_x-1}(\log x-c_x-1)-1 \lesssim \: &{\lambda}_x \: \lesssim \frac{k_x}{k_x-1}(\log x-c_x-1).\\
-\frac{k_x}{k_x-1}(\log x-c_x-1) \lesssim \: &-{\lambda}_x \: \lesssim -\frac{k_x}{k_x-1}(\log x-c_x-1)+1.\\
-(\log x-c_x-1) \lesssim \: &-{\lambda}_x \frac{k_x-1}{k_x}\: \lesssim -(\log x-c_x-1)+\frac{k_x-1}{k_x}.\\
0\lesssim \: &\log x-c_x-1-{\lambda}_x \frac{k_x-1}{k_x} \: \lesssim \frac{k_x-1}{k_x}.\\
0\lesssim \: &\frac{x}{\pi(x)}-1-{\lambda}_x \frac{k_x-1}{k_x} \: \lesssim 1.\\
1+{\lambda}_x \frac{k_x-1}{k_x} \lesssim \: &\frac{x}{\pi(x)}\: \lesssim 1+1+{\lambda}_x \frac{k_x-1}{k_x}.\\
\frac{k_xx}{k_x({\lambda}_x+2)-{\lambda}_x} \lesssim \: &\pi(x)\: \lesssim \frac{k_xx}{k_x({\lambda}_x+1)-{\lambda}_x}.\\
\frac{k_xx}{k_x({\lambda}_x+2)} \lesssim \: &\pi(x)\: \lesssim \frac{k_xx}{k_x({\lambda}_x+1)-{\lambda}_x}.\\
\frac{x}{{\lambda}_x+2} \lesssim \: &\pi(x)\: \lesssim \frac{k_xx}{{\lambda}_x(k_x-1)+k_x}.
\end{align}
\end{proof}
\subsection{A suitable variant of the prime number theorem}
\begin{theorem}
\[\pi(x)\sim \frac{x}{{\lambda}_x}\]
\end{theorem}
\begin{proof}
In $(3.1)$, we defined $k_x$ as 
\begin{align}
k_x=\pi(x)/{\pi}_2(x).
\end{align}
By $(2.7)$, we have
\begin{align}
{\lambda}_x \approx \frac{x-\pi(x)-{\pi}_2(x)}{\pi(x)-{\pi}_2(x)}.
\end{align}
It is clear from $(3.28)$ and $(3.29)$ that
\begin{align}
\text{when\: } x\to \infty, \text{ then } k_x\to \infty \text{ and } {\lambda}_x \to \infty.
\end{align}
As a result of Lemma $3.2$, we have
\begin{align}
\frac{x}{{\lambda}_x+2} \lesssim \: &\pi(x)\: \lesssim \frac{k_xx}{{\lambda}_x(k_x-1)+k_x}.
\end{align}
Dividing $(3.31)$ by $x/\log x$, then taking the limit of the right side shows that
\begin{align}
\lim_{x \to \infty} \frac{\frac{k_xx}{{\lambda}_x(k_x-1)+k_x}}{\frac{x}{{\lambda}_x}} &=\lim_{x \to \infty} \frac{k_xx{\lambda}_x}{x({\lambda}_x(k_x-1)+k_x)}.\\
&=\lim_{x \to \infty} \frac{k_x {\lambda}_x}{{\lambda}_x(k_x-1)+k_x}.\\
&=\lim_{x \to \infty} \frac{k_x}{k_x-1}.\\
&=1.
\end{align}
Similarly, taking the limit of the left side of $(3.31)$ after dividing by $x/\log x$ yields
\begin{align}
\lim_{x \to \infty} \frac{\frac{x}{{\lambda}_x+2}}{\frac{x}{{\lambda}_x}} &=\lim_{x \to \infty} \frac{x{\lambda}_x}{x({\lambda}_x+2)}.\\
&=\lim_{x \to \infty} \frac{{\lambda}_x}{{\lambda}_x+2}.\\
&=1.
\end{align}
The inequalities in $(3.31)$ combined with $(3.35)$ and $(3.38)$ imply that 
\begin{align}
\lim_{x \to \infty} \frac{\pi(x)}{\frac{x}{{\lambda}_x}} &=1.\\
\pi(x) &\sim \frac{x}{{\lambda}_x}.
\end{align}
\end{proof}
\begin{theorem}
\[{\lambda}_x \sim \log x\]
\end{theorem}
\begin{proof}
By the prime number theorem, we have
\begin{align}
\pi(x) \sim \frac{x}{\log x}.
\end{align}
It is known from Theorem $3.3$ that
\begin{align}
\pi(x) \sim \frac{x}{{\lambda}_x}.
\end{align}
Statements $(3.41)$ and $(3.42)$ imply that
\begin{align}
{\lambda}_x \sim \log x.
\end{align}
\end{proof}
\section{An asymptotic bound of ${\pi}_2(x)$}
Since ${\lambda}_x \sim \log x$ is established, we proceed to prove the following. 
\begin{theorem}
\[{\pi}_2(x) \sim (2\pm o(\log x))\frac{x}{{\log}^2 x}\]
\end{theorem}
\begin{proof}
By $(2.8)$,
\begin{align}
{\pi}_2(x)&\approx \frac{\pi(x)\big({\lambda}_x+1\big)-x}{{\lambda}_x-1}.\\
{\pi}_2(x)&\approx \frac{\pi(x)\big({\lambda}_x+1\big)}{{\lambda}_x-1}-\frac{x}{{\lambda}_x-1}.
\end{align}
Recall that in $(3.2)$, we set
\begin{align}
\pi(x)=\frac{x}{\log x-c_x},
\end{align}
where 
\begin{align}
c_x \to 1 \text{ when } x \to \infty.
\end{align}
 With $(4.2)$ and $(4.3)$, we obtain
\begin{align}
{\pi}_2(x)&\approx \frac{x({\lambda}_x+1)}{(\log x-c_x)({\lambda}_x-1)}-\frac{x}{{\lambda}_x-1}.\\
{\pi}_2(x)&\approx \frac{x}{{\lambda}_x-1}\left(\frac{{\lambda}_x+1}{\log x-c_x}-1\right).\\
{\pi}_2(x)&\approx \frac{x}{{\lambda}_x-1}\left(1+\frac{{\lambda}_x-\log x+c_x+1}{\log x-c_x}-1\right).
\end{align}
By Theorem $3.4$, 
\begin{align}
 {\lambda}_x \sim \log x.
\end{align}
This also means that 
\begin{align}
{\lambda}_x=(\log x)(1\pm o(1)).
\end{align}
Substituting ${\lambda}_x$ with $(\log x)(1\pm o(1))$ in $(4.7)$ yields
\begin{align}
{\pi}_2(x)&\approx \frac{x}{(\log x)(1\pm o(1))-1}\left(1+\frac{c_x+1\pm o(\log x)}{\log x-c_x}-1\right).\\
{\pi}_2(x)&\approx \frac{x}{(\log x)(1\pm o(1))-1}\left(\frac{c_x+1\pm o(\log x)}{\log x-c_x}\right).
\end{align}
Set
\begin{align}
C_{\pi}(x)=c_x+1\pm o(\log x),
\end{align}
then 
\begin{align}
{\pi}_2(x) \approx \frac{x}{(\log x)(1\pm o(1))-1}\left(\frac{C_{\pi}(x)}{\log x-c_x}\right).
\end{align}
Because ${\pi}_2(x)>0$, we have $C_{\pi}(x)>0$. \\
The last sentence combined with $(4.13)$ imply that
\begin{align}
{\pi}_2(x) \sim  C_{\pi}(x)\frac{x}{{\log}^2 x}.
\end{align}
Using statement $(4.4)$ to replace $c_x$ with $1$ in $(4.12)$, gives a leading constant of
\begin{align}
C_{\pi}(x)=2\pm o(\log x).
\end{align}
 With $(4.14)$ and $(4.15)$, we obtain
\begin{align}
{\pi}_2(x) &\sim  (2\pm o(\log x))\frac{x}{{\log}^2 x}.
\end{align}
\end{proof}
\section{Remarks}
\textbf{1)}  The $\lq\lq \pm "$ sign in the leading constant in $(4.16)$, is believed to be a negative $\lq\lq-"$ operator. This means that $(4.9)$ should be
\begin{align}
\text{Claim: \: \: } {\lambda}_x=(\log x)(1- o(1)).
\end{align}
To validate claim $(5.1)$, one must show that 
\begin{align}
{\lambda}_x<\log x.
\end{align}
Unfortunately, we could not provide a proof of $(5.2)$ to justify Claim $(5.1)$. For this reason, we have a $\lq\lq \pm "$ sign in our leading constant $C_{\pi}(x)$ in $(4.15)$, $(4.16)$ in lieu of a solid minus $\lq\lq - "$ sign.\\

\textbf{2)} We shall also note that, there is a possible deviation of our leading constant \\$C_{\pi}(x)=2\pm o(\log x)$, from what the elusive yet enlightening Hardy-Littlewood conjecture predicts; That is a constant value of $2C_2=1.3203.$ Nevertheless, the main result in $(4.16)$ has heretofore been out of reach.\\


\end{document}